\documentclass[12pt]{amsart}
\usepackage{fullpage}
\usepackage{amsmath}
\usepackage{algorithmicx}
\usepackage{url}
\usepackage{algpseudocode}
\newcommand{\fma}{\operatorname{fma}}
\newcommand{\ulp}{\operatorname{ulp}}
\newcommand{\abs}[1]{\left|#1\right|}
\begin{document}
\title[Finite-precision Horner approximations]{Computing accurate Horner form
approximations to special functions in finite precision arithmetic}
\author{Tor Myklebust}
\thanks{
Department of Combinatorics and Optimization, Faculty
of Mathematics, University of Waterloo, Waterloo, Ontario N2L 3G1,
Canada (e-mail: tmyklebu@csclub.uwaterloo.ca).  Research of this author
was supported by an NSERC Doctoral Scholarship.}
\date{\today}
\begin{abstract}
  In various applications, computers are required to compute approximations to
  univariate elementary and special functions such as $\exp$ and $\arctan$ to
  modest accuracy.  This paper proposes a new heuristic for automating the
  design of such implementations.  This heuristic takes a certain restricted
  specification of program structure and the desired error properties as input
  and takes explicit account of roundoff error during evaluation.
\end{abstract}

\maketitle

\begin{section}{Introduction}
  Various programming language standards, such as the Java Language
  Specification \cite{JavaLangSpec} require that a certain set of elementary
  and special functions (in this paper simply ``mathematical functions'') be
  provided for their native floating-point types and that their implementations
  are \emph{faithfully-rounded}---that is, that they produce one of the two
  machine-representable numbers bracketing the exact, real-number value.

  Projects such as LIBULTIM \cite{LibUltim} and CRlibm \cite{CRlibm} provide a
  wide array of \emph{correctly-rounded} functions---that is, they always
  produce the closest machine-representable number to the exact function value.
  CRlibm provides correctly-rounded functions in all four rounding modes
  specified by the IEEE 754 floating-point arithmetic standard.  The fdlibm
  library, developed by Sun, is a portable, widely-used math library that
  delivers faithfully-rounded results for a number of important functions.

  LIBULTIM and CRlibm make extensive use of both of table lookups and of
  conditional branches in their function implementations.  For most functions,
  the fdlibm library splits the function's domain into several intervals and
  uses a different approximation on each interval.
  
  If one wishes to compute the same mathematical function on several different
  inputs at once, it is natural to consider using the fast vector units that
  are widespread in modern computers.  It is often difficult to achieve a
  significant speedup using a vector unit with code that is rich in table
  lookups and conditional branches since different entries of the vector can
  result in different table accesses or different code paths.
  In several models of processor power consumption, using a vector instruction
  in the place of a sequence of scalar instructions can also confer a
  significant power savings \cite{Shao,Horowitz}.

  Furthermore, mathematical function evaluation is responsible for a
  substantial fraction of execution time and power use in some applications.
  \cite{Piparo} gives a brief computational study of two applications in
  high-energy physics where faster but less accurate mathematical functions
  lead to physically acceptable results within a substantially shorter
  timeframe.

  Recent projects such as SLEEF \cite{SLEEF} and Yeppp \cite{Yeppp} provide
  general-purpose, fast, and reasonably accurate mathematical function
  implementations.
  These libraries allow a programmer to compute the same function applied to a
  vector of floating-point numbers considerably faster than by applying a
  traditional implementation to each element in sequence.
  The mathematical function implementations in these libraries consist of a
  simple argument reduction that is easy to vectorise, then a polynomial
  evaluation, then a reconstruction step that is also easy to vectorise.  Both
  of these libraries focus on delivering results very quickly at the expense of
  a significantly weaker guarantee on the error between the delivered result
  and the mathematical result.

  The work in this paper is motivated by a desire to improve the error bounds
  achievable within the algorithmic framework represented by SLEEF and Yeppp
  without sacrificing the speed of the resulting implementation.

  Fundamental work on the problem of computing good machine approximations to
  functions has been done previously.  I mention only the following two recent
  papers.

  Brisebarre, Muller, and Tisserand \cite{BrisebarreMullerTisserand} give an
  algorithm based on enumerating lattice points inside a polyhedron for finding
  the polynomial with machine-representable coefficients that best
  approximates, when evaluated in real arithmetic, a given function.

  Brisebarre and Chevillard \cite{BrisebarreChevillard} give a heuristic based
  on lattice basis reduction for computing a polynomial with
  machine-representable coefficients that approximates a given function well
  when evaluated in real arithmetic.
\end{section}

\begin{section}{Overview}
  This paper gives a heuristic that tries to find polynomials with
  machine-representable coefficients that approximate a given function well
  when evaluated in machine arithmetic.  It takes as input a
  partially-specified straight-line program of fused multiply-adds (in a
  certain Horner-like form), an interval of interest, and a function that
  gives, for each machine-representable number in said interval, the range of
  acceptable function values.  The heuristic either fails or produces as output
  a fully-specified straight-line program of fused multiply-adds that produces
  an acceptable function value for every machine-representable number in the
  specified interval.

  The heuristic presented in this paper can be understood as a refinement of
  Brisebarre, Muller, and Tisserand's polyhedral approach.
  Importantly, this heuristic accounts for roundoff error that occurs
  while evaluating the function.  Further, this heuristic can take explicit
  account of any argument-reduction and reconstruction steps necessary in the
  mathematical function implementation.
  The key recent advance in optimisation technology that makes this heuristic
  practical is the exact linear optimisation package \texttt{QSopt\_ex} of
  Applegate, Cook, Dash, and Espinoza \cite{QSoptex}.

  An implementation of the heuristic in this paper, together with a
  distribution of \texttt{QSopt\_ex}, is available at
  \url{http://github.com/tmyklebu/funapprox}.
\end{section}

\begin{section}{Notation}
  C99's hexfloat notation is used extensively in this paper.
  In this notation, a finite, normal \texttt{binary32} floating-point constant
  is written unambiguously in the form \texttt{0x1.mmmmmmp+eef}, where
  \texttt{0x} is literal, each \texttt{m} is a hexadecimal digit in the
  significand, \texttt{p} is literal, \texttt{+ee} is the exponent, and
  \texttt{f} is a literal suffix indicating that the number is
  \texttt{binary32}.  See the recent C11 \cite{C11} standard's description of
  ``hexadecimal floating-point constants'' for a formal definition.
  As an example, \texttt{0x1.1f2p+1f} is the constant
  $1 \cdot 2^{1} + 1 \cdot 2^{-3} + 15 \cdot 2^{-7} + 2 \cdot 2^{-11}$.

  ``Ulp'' is a shorthand for ``unit in the last place.''
  Ignoring special cases, this is the difference between one floating-point
  number and the smallest floating-point number larger than it in magnitude.
  If $x$ is a real number, then an ``ulp of $x$'', written $\ulp x$, is the
  difference between the largest floating-point number less-than-or-equal-to
  $x$ and the smallest floating-point number larger than $x$.
  \footnote{There are several other definitions of $\ulp x$ that differ when
  $x$ is equal to or near a signed power of two.  See \cite{MullerUlp} for an
  extended discussion of various definitions of ``ulp.'' Any reasonable
  definition will do for the purpose of understanding this paper.}
\end{section}

\begin{section}{A running example}
  The following partially-specified C program will be used as a running
  example:
  \begin{verbatim}
    float sin_poly(float a) {
      float  s = a * a;
      float r5 = fmaf(s, c9, c7);
      float r4 = fmaf(s, r5, c5);
      float r3 = fmaf(s, r4, c3);
      float r2 = s * r3;
      float r1 = fmaf(a, r2, a);
      return r1;
    }
  \end{verbatim}
  Here, \texttt{fmaf} is the ``fused multiply-add'' for \texttt{binary32}
  numbers; \texttt{fmaf(a, b, c)} is the closest machine-representable number
  to $ab+c$.

  The program \texttt{sin\_poly} above attempts to compute the Horner form
  \begin{equation}\label{sinhorner}
    a + a \cdot a^2 (c_3 + a^2(c_5 + a^2(c_7 + a^2 c_9)))
  \end{equation}
  in \texttt{binary32} arithmetic.  (Note that simple multiplications such as
  \verb+s = a * a;+ may be computed as \verb+s = fmaf(a, a, 0.0f);+ if only for
  the sake of uniformity.)

  It is desired that coefficients $c_3$, $c_5$, $c_7$, and $c_9$, each in
  $[-1,1]$, are computed such that, for every \texttt{binary32} number $a$
  between $-\pi/4$ and $\pi/4$, the value returned from the call
  \texttt{sin\_poly($a$)} is within 0.65 ulp of $\sin a$.  The choice of
  $[-1,1]$ for every coefficient is arbitrary, but finite and reasonably small
  bounds on every coefficient are necessary so that the error bounds in the
  next section can give useful results.

  Two obstacles present themselves.  First, each coefficient
  ($c_3$, $c_5$, $c_7$, and $c_9$) must be a \texttt{binary32} number.
  Second, $s$, $r_5$, $r_4$, $r_3$, $r_2$, and $r_1$ are evaluated in
  \texttt{binary32} arithmetic---it is not enough that \eqref{sinhorner} is
  within 0.65 ulp of $\sin a$ for every $a \in [-\pi/4, \pi/4]$.
\end{section}

\begin{section}{Error bounds}
  The fused multiply-add (FMA) is available on many modern processors.  The
  fused multiply-add computes $ab+c$, for machine numbers $a$, $b$, and $c$,
  with only a single rounding at the end.  This calculation is written as
  $\fma(a,b,c)$.  The lack of rounding of the intermediate product $ab$ often
  results in greater accuracy in a larger computation that makes use of the
  FMA.  The IEEE standard for floating-point arithmetic \cite{ieee754},
  as of 2008, requires that $\fma(a,b,c)$ delivers the closest
  machine-representable number to $ab+c$.

  A very useful property of the FMA is that, as a univariate function of $a$,
  $\fma(a,b,c)$ is a monotone (increasing or decreasing) function.  Thus,
  given lower and upper bounds, say $\underline{f}$ and $\overline{f}$, on
  $\fma(a,b,c)$, and values of $b$ and $c$, one can use binary search to find,
  quickly and exactly, the interval of representable numbers $a$ such that
  $\underline{f} \leq \fma(a,b,c) \leq \overline{f}$.

  As mentioned, if $a$, $b$, and $c$ are machine-representable numbers, the
  fused multiply-add $\fma(a, b, c)$ is the closest machine-representable
  number to $ab+c$.  Thus, absent exponent overflow or underflow, one can bound
  \begin{equation}\label{basicbound}
    \abs{\fma(a,b,c) - (ab+c)} \leq \frac12 \ulp(ab+c).
  \end{equation}

  If $a$ and $c$ are machine approximations to real numbers
  $a + \delta a$ and $c + \delta c$ but $b$ is a machine-representable number
  known exactly, one can use the triangle inequality to bound
  \begin{equation}\label{fmaerr}
    \begin{array}{rl}
         & \abs{\fma(a,b,c) - ((a+\delta a)b+(c + \delta c))}\\
    \leq & \frac12 \ulp(ab+c) + \abs{\delta c + b \delta a}\\
    \leq & \frac12 \ulp(ab+c) + \abs{\delta c} + \abs{b \delta a}.
    \end{array}
  \end{equation}
  Thus, if one has bounds on $ab+c$, $\abs{\delta a}$, and $\abs{\delta c}$,
  one can compute an explicit bound on the difference between $\fma(a,b,c)$ and
  the exact-arithmetic result $(a+\delta a)b + (c+\delta c)$.

  Applying this bound naively (i.e. taking no account of possible cancellation)
  to the example of \texttt{sin\_poly(0.5f)}, one
  can compute $s = 0.25$ and then find
  $$
    r_5 + \delta r_5 = 0.25 c_9 + c_7
  $$
  with
  $\abs{\delta r_5} \leq \frac12 \ulp{1.25} = 2^{-25}$ and
  $\abs{r_5} \leq 1.25 + 2^{-25}$.

  Then
  $$
     r_4 + \delta r_4
   = c_5 + 0.25 c_7 + 0.25^2 c_9
  $$
  with
  $$   \abs{\delta r_4}
  \leq \frac12 \ulp(1.3125 + 2^{-25}) + \abs{0.25 \delta r_5}
  \leq 2^{-25} + 2^{-27}
  $$
  and
  $\abs{r_4} \leq 1.3125 + 2^{-25} + 2^{-27}$.

  Bounds on $\abs{\delta r_3}$, $\abs{r_3}$, $\abs{\delta r_2}$,
  $\abs{r_2}$, $\abs{\delta r_1}$, and $\abs{r_1}$ can be computed analogously.

  Note also that, given $a$, bounds that must be satisfied by $r_1$ can be
  computed.
  For $a = \texttt{0.5f}$, we desire that $r_1$ is within $0.65$ ulp of
  $\sin \frac12$.  Since $\sin \frac12$ is roughly
  $\texttt{0x1.eaee8744b0p-2f}$, this means that we desire
  $$
    \texttt{0x1.eaee86p-2f} \leq r_1 \leq \texttt{0x1.eaee88p-2f}.
  $$
  Note that these lower and upper bounds are adjacent \texttt{binary32} numbers.
  In some other cases, such as $a = \texttt{0.625f}$, these lower and upper
  bounds are equal.

  Since $r_1$ is computed using only $r_2$ and $a$, this implies bounds on
  $r_2$.  One can use
  binary search to find the smallest and largest \texttt{binary32} numbers
  such that
  $$
    \texttt{0x1.eaee86p-2f} \leq \fma(a, r_2, a) \leq \texttt{0x1.eaee88p-2f}.
  $$
  This yields the bounds
  $$
    \texttt{-0x1.5117aep-5f} \leq r_2 \leq \texttt{-0x1.511790p-5f}.
  $$
  These bounds are nonadjacent \texttt{binary32} numbers.
  Since $r_2$ is computed using only $r_3$ and $s$, this implies bounds on
  $r_3$:
  $$
    \texttt{-0x1.5117aep-3f} \leq r_3 \leq \texttt{-0x1.511790p-3f}.
  $$
  Unfortunately, $r_3$ is computed using the coefficient $c_3$, so similar
  bounds on $r_4$ cannot be obtained.  However, if $c_3$ is fixed, this can
  be done.
\end{section}

\begin{section}{Formulating linear constraints}
  Given a value of the abscissa $a$, the value of $s$ can be computed
  directly and acceptable bounds on $r_1$ can be derived from the value of
  $\sin a$.
  As in the previous section, upper and lower bounds on the difference
  between the computed value \texttt{sin\_poly(a)} and the exact-arithmetic
  value of the Horner form \eqref{sinhorner} can be found.

  Suppose this difference is at least $\underline{\delta}$ and at most
  $\overline{\delta}$.  Then the coefficients $c_{3, \dots, 9}$ must satisfy
  the \emph{linear} inequalities
  \begin{equation}\label{diffs}
    \begin{array}{rl}
         & \sin(a) - 0.65 \ulp(\sin(a)) + \underline{\delta}\\
    \leq & a + a\cdot s(c_3 + s(c_5 + s(c_7 + sc_9)))\\
    \leq & \sin(a) + 0.65 \ulp(\sin(a)) + \overline{\delta}.
    \end{array}
  \end{equation}
  An acceptable list of coefficients $c_{3, \dots, 9}$ must satisfy
  \eqref{diffs} for every $a$ of interest.  However, one need not formulate all
  such constraints in the beginning.  One can generate the inequalities
  \eqref{diffs} for some small subset of the points of interest, find a
  solution, and then look for points $a$ not yet considered for which
  \eqref{diffs} is violated.  This is often called a cutting-plane method.

  Indeed, a classical theorem of Helly \cite{Helly} (see also \cite{Schrijver})
  implies that a (possibly very large) system of linear inequalities in $n$
  variables is either satisfiable or there exists an unsatisfiable subsystem of
  at most $n+1$ inequalities.

  Given a list of linear inequalities that must be satisfied by some variables
  that take on real values, one can compute lower and upper bounds on the
  variables using linear optimisation.

  The linear optimisation problems that arise from \eqref{diffs} are very
  ill-conditioned.  The left-hand side is a Vandermonde matrix and the bounds
  on each row of this Vandermonde, called $\underline{\delta}$ and
  $\overline{\delta}$ in \eqref{diffs}, are often very close together.
  Conventional inexact linear optimisation, as implemented in systems such as
  Gurobi \cite{Gurobi} and IBM's CPLEX \cite{cplex}, can give meaningfully
  incorrect bounds and even confuse feasible systems of linear inequalities
  with infeasible systems.  Therefore, fast, exact linear optimisation is
  necessary.

  Note also that, if $c_3$ is fixed to some value, different (likely tighter)
  linear inequalities on $c_5$, $c_7$, and $c_9$ may be formulated since, for
  each abscissa $a$, the interval of acceptable values of $r_4$ can now be
  derived.
\end{section}

\begin{section}{The heuristic}
  The state of the heuristic has several components:
  \begin{itemize}
    \item A nonempty list of \emph{test points}.
    \item Finite lower and upper bounds on each variable.
  \end{itemize}

  The heuristic runs the following loop until something fails:
  \begin{enumerate}
    \item Find a solution $c$ such that $p_c(x)$ is an acceptable rounding of
      $f(x)$ for every test point $x$.
    \item Find a point $x$ such that $p_c(x)$ is not an acceptable rounding of
      $f(x)$ and add $x$ to the list of test points.
    \item Go to 1.
  \end{enumerate}

  Step 1 is done roughly as in Figure \ref{fig:step1}, but with some refinements
  described later.
  If step 1 fails, the heuristic reports failure.
  This does not necessarily imply that the problem is infeasible.

  Step 2 is done by trying every abscissa of interest until at least one
  results in an unacceptable function value.  This is only practical for
  smaller domains, such as those arising from IEEE \texttt{binary32} functions.

  \begin{figure}
  \begin{algorithmic}
    \State Fix an ordering of the coefficients $[c_1, \dots, c_k]$.
    \Loop{\,\,some number of times}
      \For{$i = 1$ to $k$}
        \State If $i = k$, report success.
        \State Compute lower and upper bounds on $c_i$ by exact
          linear optimisation---say
          $\underline{c_i} \leq c_i \leq \overline{c_i}$.
        \State If infeasible, break.
        \State Fix $c_i$ to a randomly-chosen representable number between
          $\underline{c_i}$ and $\overline{c_i}$.
      \EndFor
    \EndLoop
    \State Fail.
  \end{algorithmic}
  \caption{Pseudocode for Step 1.}\label{fig:step1}
  \end{figure}

  Choosing a distribution other than the uniform distribution on
  $[\underline{c_i}, \overline{c_i}]$ may result in better performance.  I use
  the average of two uniform samples on $[\underline{c_i}, \overline{c_i}]$;
  better choices may exist.

  It is also wasteful to take only a single sample of $c_{i+1}, \dots, c_k$
  after fixing $c_i$.  The linear optimisation problems later in the
  \textbf{for} loop tend to solve faster than those earlier in the loop
  Further, a single bad coefficient choice later in the variable fixing process
  can scuttle a good choice of an earlier coefficient.  Instead, I recursively
  try a constant number (four) of choices of $c_{i+1}$ after fixing $c_i$.  The
  first choice of $c_i$ made is always the value of $c_i$ that most recently
  yielded an acceptable list of coefficients.
\end{section}

\begin{section}{Examples}
  This section gives examples of C functions that compute $\arctan$ on
  $[-1, 1]$ and on $(-\infty, \infty)$ using fused multiply-add.  The
  $\arctan$ function was chosen because it admits an especially simple
  argument reduction to $[-1, 1]$, yet some care must be taken to
  get a faithfully-rounded result on that interval.

  \begin{subsection}{Arctangent on $[-1,1]$}
    The C function in Figure \ref{fig:atan}
    is a faithfully-rounded approximation to $\arctan
    x$ for $x$ a \texttt{binary32} number in $[-1, 1]$.
    \begin{figure}
    \begin{verbatim}
      float atan_poly(float a) {
        float s = a * a;
        float r = 0x1.6d2026p-9f;
        r = fmaf(r, s, -0x1.03f2d4p-6f);
        r = fmaf(r, s,  0x1.5beeb4p-5f);
        r = fmaf(r, s, -0x1.33194ep-4f);
        r = fmaf(r, s,  0x1.b403a8p-4f);
        r = fmaf(r, s, -0x1.22f5c2p-3f);
        r = fmaf(r, s,  0x1.997748p-3f);
        r = fmaf(r, s, -0x1.5554d8p-2f);
        r = r * s;
        return fmaf(r, a, a);
      }
    \end{verbatim}
    \caption{A faithful approximation to $\arctan$ on $[-1,1]$.}
    \label{fig:atan}
    \end{figure}
    On every \texttt{binary32} number $a$ in $[-1, 1]$,
    $\mathtt{atan\_poly}(a)$ produces a result that differs from $\arctan a$
    by less than 0.95 ulp.  This was computed in about five minutes from the
    partial straight-line program in Figure \ref{fig:partatan}, fixing
    variables in the order $c_3, c_5, c_7, c_9, c_{11}, c_{13}, c_{15},
    c_{17}$.
    \begin{figure}
    \begin{verbatim}
    float atan_poly(float a) {
      float s = a * a;
      float r = c17;
      r = fmaf(r, s, c15);
      r = fmaf(r, s, c13);
      r = fmaf(r, s, c11);
      r = fmaf(r, s, c9);
      r = fmaf(r, s, c7);
      r = fmaf(r, s, c5);
      r = fmaf(r, s, c3);
      r = r * s;
      return fmaf(r, a, a);
    }
    \end{verbatim}
    \caption{A partially-specified approximation to $\arctan$.}
    \label{fig:partatan}
    \end{figure}

    By way of comparison, SLEEF's implementation \texttt{xatanf} has a maximum
    error of roughly 1.773 ulp on $[-1, 1]$, which drops to roughly 1.707 ulp
    if \texttt{mlaf} is replaced by \texttt{fmaf}.

    The Sollya tool \cite{Sollya}, when given input
    \begin{verbatim}
    fpminimax(atan(x), [|1,3,5,7,9,11,13,15,17|],
        [|24,24,24,24,24,24,24,24,24|], [1e-6, 1]);
    \end{verbatim}
    produces the coefficients given in Figure \ref{fig:sollya}.
    \begin{figure}
    \begin{verbatim}
    float c17 = 2.90188402868807315826416015625e-3f;
    float c15 = -1.62907354533672332763671875e-2f;
    float c13 = 4.3082617223262786865234375e-2f;
    float c11 = -7.5408883392810821533203125e-2f;
    float c9 = 0.1066047251224517822265625f;
    float c7 = -0.14209578931331634521484375f;
    float c5 = 0.19993579387664794921875f;
    float c3 = -0.3333314359188079833984375f;
    \end{verbatim}
    \caption{Sollya's coefficients for approximating $\arctan$ on $[-1,1]$.}
    \label{fig:sollya}
    \end{figure}
    Sollya's \texttt{fpminimax} uses an implementation of the method of
    Brisebarre and Chevillard.  When the coefficients of this polynomial are
    substituted into Figure \ref{fig:partatan}, the resulting function has a
    maximum error of roughly 1.067 ulp.  The proposed heuristic therefore gives
    an improvement in the maximum error of about 0.117 ulp over Sollya.
  \end{subsection}

  \begin{subsection}{Arctangent everywhere}
    The proposed approach allows one to specify the function to be approximated
    by means of a function taking an argument $x$ and returning the range of
    acceptable values of $f(x)$.  This can be used to find approximations that
    make a given mathematical function implementation (including, say, an
    argument reduction step) more accurate.  Consider the partially-specified
    implementation of $\arctan$ (due to N. Juffa \cite{JuffaQuestion}) in
    Figure \ref{fig:juffaatan}.

    \begin{figure}
    \begin{verbatim}
    float juffa_atanf(float a) {
        float r, t;
        t = fabsf(a);
        r = t;
        if (t > 1.0f) r = 1.0f / r;
        r = atan_poly(r);
        if (t > 1.0f) r = fmaf(0x1.ddcb02p-1f, 0x1.aee9d6p+0f, -r);
        r = copysignf(r, a);
        return r;
    }
    \end{verbatim}
    \caption{N. Juffa's $\arctan$ skeleton.  This calls \texttt{atan\_poly}
      from Figure \ref{fig:partatan}}\label{fig:juffaatan}
    \end{figure}

    Suppose one desires a \texttt{binary32} approximation of $\arctan$ that is
    within $1.2$ ulp of $\arctan$ everywhere.  This can be cast as an
    approximation problem on $[0, 1]$.  For every \texttt{binary32} number
    $x$ in $[0, 1]$, one wants
    \begin{itemize}
      \item the straight-line program given by \texttt{atanf\_poly} in
        \ref{fig:juffaatan} to yield a result within $1.2$ ulp of $\arctan(x)$,
        and
      \item for every $y \in (1, \infty)$ such that $1/y$ rounds to $x$,
        for the straight-line program given by \texttt{atanf\_poly} in
        \ref{fig:juffaatan} to yield a result such that, after line 18 is run,
        \texttt{r} is within $1.2$ ulp of $\arctan(y)$.
    \end{itemize}
    One can find the interval of machine-representable $y \in (1, \infty)$ such
    that $1/y$ rounds to $x$ by binary search.
    It is straightforward to write a function that computes, given \texttt{x},
    the range of acceptable values of \texttt{atanf\_poly(x)} under the above
    two conditions.

    \begin{figure}
    \begin{verbatim}
    float c17 = 0x1.686c56p-9;
    float c15 = -0x1.01dec8p-6;
    float c13 = 0x1.5a901p-5;
    float c11 = -0x1.32b648p-4;
    float c9 = 0x1.b3f558p-4;
    float c7 = -0x1.22f90cp-3;
    float c5 = 0x1.99782cp-3;
    float c3 = -0x1.5554d8p-2;
    \end{verbatim}
    \caption{Coefficients for Figure \ref{fig:juffaatan}.}
    \label{fig:atancoeffs}
    \end{figure}
    The proposed heuristic finds the coefficients in Figure
    \ref{fig:atancoeffs} in about ten minutes.  These coefficients yield an
    error less than 1.1978 ulp everywhere.  By way of comparison, Sollya's
    polynomial yields a maximum error of more than 1.535 ulp.  I also tried
    asking for a maximum error of 1.1 ulp, but the heuristic failed to find a
    solution.
  \end{subsection}
\end{section}

\begin{section}{Concluding remarks}
  This paper presented a heuristic for designing floating-point approximations
  to univariate elementary functions.  It takes as input a straight-line
  program structure, an interval, and a function giving the range of acceptable
  values for each input in the interval.
  If the heuristic succeeds, it outputs a straight-line program implementing
  the function on the desired interval that satisfies the desired error bound.

  This paper also presented two nontrivial examples of approximations to
  $\arctan$.  These produced \texttt{binary32} numbers approximating the
  arctangent of \texttt{binary32} numbers using \texttt{binary32} arithmetic.

  For larger floating-point types such as IEEE \texttt{binary64}, it is not
  clear to me how to prove in a reasonable amount of time that a given
  straight-line program either satisfies the desired error bound on an interval
  or to find an abscissa where it fails.  Moreover, this approach might not
  scale to the higher-degree polynomials necessary for such approximations.  I
  leave these considerations to future work.
\end{section}

\begin{section}{Acknowledgements}
  I thank Norbert Juffa, Levent Tun\c{c}el, and Jennifer Wong for their
  excellent feedback on early drafts of this paper.  I thank Norbert Juffa
  also for posing this problem to me and for sharing with me some of his
  considerable knowledge and wisdom on this subject. 
\end{section}

\addcontentsline{toc}{chapter}{References}
\bibliographystyle{plain}
\bibliography{funapprox}
\nocite{*}
\end{document}